\documentclass[12pt,twoside,a4paper]{article}
\usepackage{t1enc} 
\usepackage{times}
\usepackage[dvips]{epsfig}
\usepackage{amssymb}
\usepackage{amsmath}
\usepackage{longtable}
\usepackage{makeidx}
\usepackage{euscript}
\usepackage{enumerate}
\usepackage{fancyhdr}
\usepackage{hyperref}
\DeclareMathAlphabet\scr{U}{scr}{m}{n}
\SetMathAlphabet\scr{bold}{U}{scr}{b}{n}
  \DeclareFontFamily{U}{scr}{\skewchar\font'177}%
  \DeclareFontShape{U}{scr}{m}{n}{<-6>rsfs5<6-8>rsfs7<8->rsfs10}{}%
  \DeclareFontShape{U}{scr}{b}{n}{<-6>rsfs5<6-8>rsfs7<8->rsfs10}{}%

\newcommand{\ev}{\mathbb E}

\newcommand{\pp}{\mathbb P}

\newtheorem{satz}{Theorem}[section]
\newtheorem{prop}[satz]{Proposition}
\newtheorem{theorem}[satz]{Theorem}

\newtheorem{lemma}[satz]{Lemma}

\newtheorem{assumption}[satz]{Assumption}
\newtheorem{@definition}[satz]{Definition}

\newtheorem{@aufgabe}{Aufgabe}

\newtheorem{@bsp}[satz]{Beispiel}

\newenvironment{mat}{\left (\begin{matrix} } {\end{matrix}\right ) }

\newcommand{\bmat}{\begin{mat}}
\newcommand{\emat}{\end{mat}}
\newcommand{\be}{\begin{enumerate}}
\newcommand{\ee}{\end{enumerate}}
\newcommand{\beq}{\begin{equation}}
\newcommand{\eeq}{\end{equation}}
\newcommand{\bea}{\begin{eqnarray}}
\newcommand{\eea}{\end{eqnarray}}
\newcommand{\beaa}{\begin{eqnarray*}}
\newcommand{\eeaa}{\end{eqnarray*}}

\newcommand{\ep}{\hfill $\Box$}

\renewcommand{\epsilon}{\varepsilon}
\renewcommand{\phi}{\varphi}
\renewcommand{\rho}{\varrho}

\title{Metastability in the generalized Hopfield model with finitely many patterns}
\author{Mykhaylo Shkolnikov}
\date{}

\begin{document}
\maketitle

\begin{abstract}
This paper continues the study of metastable behaviour in disordered mean field models initiated in \cite{bbi}, \cite{begk1}. We consider the generalized Hopfield model with finitely many independent patterns $\xi_1,\dots,\xi_p$ where the patterns have i.i.d. components and follow discrete distributions on $[-1,1]$. We show that metastable behaviour occurs and provide sharp asymptotics on metastable exit times and the corresponding capacities. We apply the potential theoretic approach developed by Bovier et al. in the space of appropriate order parameters and use an analysis of the discrete Laplacian to obtain lower bounds on capacities. Moreover, we include the possibility of multiple saddle points with the same value of the rate function and the case that the energy surface is degenerate around critical points. 
\end{abstract}

\section{Introduction}

We consider the generalized Hopfield model with finitely many patterns which was introduced for a special case in \cite{pf} and reintroduced for a more general case in \cite{hop}. As explained in both references it has applications in the theory of spin glasses and neural networks. In the language of mean field models the Hamiltonian in this model depends on a fixed number $p\geq1$ of random macroscopic quantities where the randomness is introduced through independent random interaction patterns $\xi_1,\dots,\xi_p$. We assume that each pattern $j$ has i.i.d. components $\xi_j(i)$, $1\leq i\leq n$ following discrete distributions on $[-1,1]$. With these patterns the Gibbs measure $\mu_n$ on spin configurations $\sigma=(\sigma(i))_{i=1}^n\in\{-1,1\}^n$ is given by
\begin{eqnarray}
\mu_n(\sigma)=\frac{1}{Z_n}e^{-\beta H_n(\sigma)} 
\end{eqnarray} 
with Hamiltonian
\begin{eqnarray}
H_n(\sigma)=-n\cdot v\left(\frac{1}{n}\left\langle \sigma,\xi_1\right\rangle,\dots,\frac{1}{n}\left\langle\sigma,\xi_p\right\rangle\right).
\end{eqnarray}
Here, $Z_n$ is the normalization constant making $\mu_n$ a probability measure, $\beta>0$ is a parameter which has the physical interpretation of the inverse temperature and $v$ is a non-negative $C^2$ function which stands for the potential of the model.\\\quad\\ 
We observe that the Hamiltonian depends only on the order parameters 
\begin{eqnarray}
X=(X_1,\dots,X_p)=\left(\frac{1}{n}\left\langle \sigma,\xi_1\right\rangle,\dots,\frac{1}{n}\left\langle \sigma,\xi_p\right\rangle\right). 
\end{eqnarray}
By Theorem V in \cite{c} the latter satisfy a large deviation principle in the limit $n\rightarrow\infty$ with sequence $n$ and rate function
\begin{eqnarray}
I=-\beta v + L^* + c
\end{eqnarray}
where $L^*$ is the Fenchel-Legendre transform of $L(t)=\ev[\log\cosh(t\cdot\xi_.(1))]$ and $c=\lim_{n\rightarrow\infty} \frac{1}{n}\log Z_n$. In other words, denoting the distribution of $(X_1,\dots,X_p)$ under $\mu_n$ by $Q_n$ we have
\begin{eqnarray}
Q_n(X)=e^{-nI(X)}\kappa_n(X)(1+o(1))
\end{eqnarray}
with a subexponential correction term $\kappa_n(X)$. We make the following assumption on $I$:
{\assumption Let $v$ be such that $I$ admits several local minima with different values of $I$ and assume that for each two minima the critical points of $I$ with the minimal value of $I$ over all paths connecting the two minima form a finite set consisting of saddle points of $I$. Moreover, let the distributions of the patterns be such that $I$ is finite on a subset of $[-1,1]^p$ with non-empty interior and $\min_{\partial\{x\in[-1,1]^p:\;I(x)<\infty\}} I$ is larger than the value of $I$ at any of the saddle points mentioned above.}\\\quad\\
{\bf Remark.} To clarify the conditions on $I$, we give examples of two situations in which the former are satisfied. Let none of the distributions of the components of the patterns be the Dirac mass at $0$. Then the set 
\begin{eqnarray*}
\{x\in[-1,1]^p:\;I(x)<\infty\}=\{x\in[-1,1]^p:\;L^*(x)<\infty\}
\end{eqnarray*}
is a convex subset of $[-1,1]^p$ with non-empty interior. Consequently, $L^*$ has a unique minimum at a point $x_0$. If $v$ is the potential of the classical Hopfield model
\begin{eqnarray*}
v(x)=x_1^2+\dots+x_p^2,
\end{eqnarray*}
then for $\beta$ exceeding the largest eigenvalue of the Hessian of $L^*$ at $x_0$, but small enough such that the condition on $\min_{\partial\{x\in[-1,1]^p:\;I(x)<\infty\}} I$ is satisfied, the rate function $I$ will have multiple local minima, all in the interior of $\{x\in[-1,1]^p:\;I(x)<\infty\}$. An other example of a situation in which $I$ has multiple local minima in the interior of $\{x\in[-1,1]^p:\;I(x)<\infty\}$ is the case of the classical Hopfield model with random fields corresponding to the potential
\begin{eqnarray*}
v(x)=\sum_{i=1}^{p-1} x_i^2 + x_p 
\end{eqnarray*}  
with same constraints on $\beta$. The existence of such a regime of $\beta$ was shown in \cite{begk1} for the case of the random field Curie-Weiss model, i.e. $p=2$ and the components of the first pattern being identically equal to $1$. Finally, we observe that we need to break the symmetry of $I$ to obtain local minima with different values of $I$ separated by saddle points of $I$. This can be done by choosing the distributions of the patterns to be non-symmetric.\\\quad\\
The generalized Hopfield model can be equipped naturally with a reversible Markovian dynamics in the following way. We let the configuration $\sigma$ evolve in disrete time with Markovian transition rates
\begin{eqnarray*}
p_n(\sigma,\sigma^i)=\frac{1}{n}\exp(-\beta(H_n(\sigma^i)-H_n(\sigma))_+),\\
p_n(\sigma,\sigma)=1-\frac{1}{n}\sum_{i=1}^n \exp(-\beta(H_n(\sigma^i)-H_n(\sigma))_+)
\end{eqnarray*}  
where $\sigma^i$ is the configuration obtained from $\sigma$ by flipping its $i$-th coordinate. A simple computation shows that the detailed balance condition holds. We observe that the described dynamics induces a corresponding (in general non-Markovian) dynamics of the vector $X$ of order parameters with transition probabilities
\begin{eqnarray}
r_n(X,X')=\frac{1}{Q_n(X)}\sum_{\sigma:\pi(\sigma)=X} \mu_n(\sigma)\sum_{\sigma':\pi(\sigma')=X'} p_n(\sigma,\sigma')
\end{eqnarray}
where 
\begin{eqnarray}
\pi:\;\{-1,1\}^n\rightarrow\Lambda_n\subset[-1,1]^p,\quad\sigma\mapsto(X_1,\dots,X_p) 
\end{eqnarray}
is the projection on the order parameters.\\\quad\\ 
Let $m=m_n$ be a local, but not a global, minimum of the rate function $I$ on $\Lambda_n$ and $M=M_n$ be the set of all local minima on $\Lambda_n$ with rate smaller than $I(m)$. Note that for large $n$ the critical points of $I$ on $\Lambda_n$ are just $\Lambda_n$-valued lattice approximations of the critical points of $I$ on $[-1,1]^p$. For this reason, with an abuse of notation we will denote both types of critical points with the same letters. We start the Markov chain in a configuration in $\pi^{-1}(m)$ and study the entrance time $\tau_{\pi^{-1}(M)}$ into the set $\pi^{-1}(M)$. Under the assumption 1.1 we compute the sharp asymptotics of such entrance times in the limit $n\rightarrow\infty$. We observe metastable behaviour in the sense that up to a correction term $\tau_{\pi^{-1}(M)}$ grows exponentially in $n$. It turns out that the metastable time scale depends on $I(z)-I(m)$ where $z$ belongs to the set ${\scr Z}$ of saddle points on paths from $m$ to $M$ for which the value of $I$ is minimal. Moreover, the correction term $c(n)$ in front of the exponential depends only on the local geometry of the energy landscape near the minimum $m$ and the saddle points in ${\scr Z}$ and the transition probabilities at these critical points.\\\quad\\
The main object for the study of entrance times will be the extended vector of order parameters $Y$ which we introduce next. Denote by $A$ the finite set of possible values of the vector $(\xi_1(1),\dots,\xi_p(1))$ and define the extended vector of order parameters by
\begin{eqnarray}
Y=(Y^{\pm}_a:\;a\in A)
\end{eqnarray}
where
\begin{eqnarray}
Y^+_{a}=\frac{1}{n}|\{i:\;(\xi_1(i),\dots,\xi_p(i))=a,\;\sigma(i)=1\}|,\quad a\in A,\\
Y^-_{a}=\frac{1}{n}|\{i:\;(\xi_1(i),\dots,\xi_p(i))=a,\;\sigma(i)=-1\}|,\quad a\in A.
\end{eqnarray}
Denote the dimension of $Y$ by $L=2|A|$, the distribution of $Y$ by $\widehat{Q}_n$ and the transition rates of the induced dynamics of $Y$ by
\begin{eqnarray}
\widehat{r}_n(Y,Y')=\frac{1}{\widehat{Q}_n(Y)}\sum_{\sigma:\pi_1(\sigma)=Y} \mu_n(\sigma)\sum_{\sigma':\pi_1(\sigma')=Y'} p_n(\sigma,\sigma').
\end{eqnarray}
In contrast to the dynamics of $X$, the dynamics of $Y$ is a reversible Markovian dynamics. This is due to the fact that for any fixed disorder the value of $Y$ is enough to determine the possible transition directions from $Y$. Indeed, in each step either $Y$ does not change at all or there is a unique $a\in A$ such that one of $Y_a^{\pm}$ increases by $\frac{1}{n}$ and one decreases by $\frac{1}{n}$ and $Y_{a'}$ stays the same for $a'\neq a$. Furthermore, for an $a\in A$ the described transition has non-zero probability if and only if the resulting vector $Y'$ is an element of $[0,1]^L$. Moreover, knowing $Y$ and $Y'$ we know if the transition corresponds to a flip from $+$ to $-$ or to a flip from $-$ to $+$ and the transition probability is given by
\begin{eqnarray*}
\widehat{r}_n(Y,Y')=\frac{1}{n}\exp(-\beta n(v(\pi_2(Y))-v(\pi_2(Y')))_+)\widehat{Q}_n(Y)^{-1}\sum_{\sigma:\pi_1(\sigma)=Y}\mu_n(\sigma)\\
\cdot|\{i:\pi_1(\sigma^i)=Y'\}|=\exp(-\beta n(v(\pi_2(Y))-v(\pi_2(Y')))_+)Y^{\pm}_a,
\end{eqnarray*}
respectively. Thus, $\widehat{r}_n(Y,Y')$ can be written as a function of only $Y$ and $Y'$, so the dynamics is Markovian.\\\quad\\ 
In addition, one observes that the projection $\pi$ decomposes in two projections
\begin{eqnarray*}
\pi:\;\{-1,1\}^n \longrightarrow \widehat{\Lambda}_n \longrightarrow \Lambda_n,\\
\sigma\;\;\stackrel{\pi_1}{\mapsto}\;\;Y\;\;\stackrel{\pi_2}{\mapsto}\;\;X.
\end{eqnarray*}
Our main result is the following 
{\theorem Under assumption 1.1 it holds 
\begin{eqnarray}
\ev^{\pi^{-1}(m)}[\tau_{\pi^{-1}(M)}]=\ev^{\pi_2^{-1}(m)}[\tau_{\pi_2^{-1}(M)}]=c(n)\exp(n(I(z)-I(m)))(1+o(1))
\end{eqnarray}
with a subexponential correction term $c(n)$ where the expectation is taken with respect to the initial distribution
\begin{eqnarray}
\nu(Y)=\frac{\widehat{Q}_n(Y)\pp^Y(\tau_{\pi_2^{-1}(M)}<\tau_{\pi_2^{-1}(m)})}
{\sum_{Y'\in \pi_2^{-1}(m)}\widehat{Q}_n(Y')\pp^{Y'}(\tau_{\pi_2^{-1}(M)}<\tau_{\pi_2^{-1}(m)})}
\end{eqnarray}
on $\pi_2^{-1}(m)$ and the corresponding distribution on $\pi^{-1}(m)$, respectively. Moreover, the correction term $c(n)$ can be computed and is given in proposition 3.1.}\\\quad\\
The statement of Theorem 1.1 is of the same type as the main result of \cite{bbi} where the correction term was computed explicitly in the case of the Curie-Weiss model with continuous random fields which is an extension of the metastability result in the Curie-Weiss model with discrete random fields gived in \cite{begk1}. We extend the latter result in a different direction treating the case of any number of discrete patterns and a general potential. In our setting the correction term has a complicated structure (see Proposition 3.1) due to the fact that several saddle points with the same value of $I$ may occur, the Hessian $H$ of $I$ may be degenerate at critical points of $I$ and the correction term in the large deviation principle for $Q_n$ is not known in general.\\\quad\\
The proof of Theorem 1.1 relies on the potential theoretic approach to metastability developed by Bovier et al. in \cite{begk1}, \cite{begk2}, \cite{begk3}. We define the Dirichlet form with respect to the weights $\widehat{Q}_n(Y)\widehat{r}_n(Y,Y')$ by 
\begin{eqnarray}
\widehat{d}(f)=\frac{1}{2}\sum_{Y,Y'} \widehat{Q}_n(Y)\widehat{r}_n(Y,Y')(f(Y')-f(Y))^2
\end{eqnarray}
and the corresponding capacity $\widehat{cap}(A,B)$ as the infimum of $d$ over all functions $f$ with $f|_A\equiv1$, $f|_B\equiv0$. We recall that the infimum is attained by the unique harmonic function $\widehat{\Phi}$ satisfying the given boundary conditions and that the following formula from potential theory holds:  
\begin{eqnarray}
\ev^{\pi_2^{-1}(m)}[\tau_{\pi_2^{-1}(M)}]=\frac{1}{\widehat{cap}(\pi_2^{-1}(m),\pi_2^{-1}(M))}\sum_{Y} \widehat{Q}_n(Y)\widehat{\Phi}(Y)
\end{eqnarray}
whereby $A=\pi_2^{-1}(m)$ and $B=\pi_2^{-1}(M)$. As usually the main part of the proof is the computation of the sharp asymptotics of 
\begin{eqnarray*}
\widehat{cap}(\pi_2^{-1}(m),\pi_2^{-1}(M)). 
\end{eqnarray*}
These are usually derived using a (simple) upper bound of the form $\widehat{d}(g)$ for a suitable function $g$ and a (hard) lower bound resulting from the Raighley's comparison principle for Dirichlet forms or the Berman-Konsowa variational principle for capacities (see \cite{bk}, \cite{bbi}). In this paper we derive the upper bound as explained using a function $g$ introduced in \cite{bbi} for the study of metastability in the random field Curie-Weiss model. One of the original contributions of this paper is the derivation of the lower bound. Here, we prove an estimate directly on $|\widehat{d}(g)-\widehat{cap}(\pi_2^{-1}(m),\pi_2^{-1}(M))|$ by analyzing the discrete Laplacian with respect to the weights $\widehat{Q}_n(Y)\widehat{r}_n(Y,Y')$. Moreover, we do not assume non-degeneracy of $H$ at critical points as was typically done in previous work on metastability.\\\quad\\
The paper is organized as follows. In section 2 we prove the upper bound on $\widehat{cap}(\pi_2^{-1}(m),\pi_2^{-1}(M))$ in Lemma 2.1 and show that it is sharp in Theorem 2.3. In section 3 we prove Theorem 1.2 and provide the exact value of the correction term $c(n)$ using the sharp asymptotics of the capacity given in Theorem 2.3.     

\section{Asymptotics of the capacity}

In this section we derive the sharp asymptotics of the capacity. This is done by obtaining an upper bound in Lemma 2.1 and showing that it is sharp in Theorem 2.3. To this end, we need a better understanding of the energy landscape for the extended vector of order parameters $Y$.\\\quad\\
{\bf 2.1. Large deviation principle for the extended vector of order parameters.} We start our analysis with the derivation of a large deviation principle for $Y$ from the large deviation principle for $X$. To this end, we observe that for a fixed realization of the patterns the value of $Y_a^+ + Y_a^-$ is fixed for every $a\in A$. Since
\begin{eqnarray}
X=\sum_{a\in A} a(Y_a^+ - Y_a^-),
\end{eqnarray}
the mapping $\pi_2$ maps at most $(n+1)^{|A|}$ points in $\widehat{\Lambda}_n$ to the same point in $\Lambda_n$. Moreover, note that $\widehat{Q}_n(Y)$ depends only on $Q_n(\pi_2(X))$. Hence, the extended order parameters satisfy a large deviation principle with sequence $n$ and rate function $\widehat{I}\equiv I\circ\pi_2$. In other words, 
\begin{eqnarray}
\widehat{Q}_n(Y)=e^{-n\widehat{I}(Y)}\widehat{\kappa}_n(Y)(1+o(1))
\end{eqnarray}  
where 
\begin{eqnarray}
\widehat{\kappa}_n(Y)=\frac{\kappa_n(\pi_2(Y))}{|\{Y':\:\pi_2(Y')=\pi_2(Y)\}|}. 
\end{eqnarray}
\quad\\
{\bf 2.2. Upper bound on the capacity.} Since the dynamics of $Y$ is a reversible Markov chain as explained in the introduction, we have
\begin{eqnarray}
\ev^{\pi^{-1}(m)}[\tau_{\pi^{-1}(M)}]=\ev^{\pi_2^{-1}(m)}[\tau_{\pi_2^{-1}(M)}]=\frac{1}{\widehat{cap}(\pi_2^{-1}(m),\pi_2^{-1}(M))}\sum_{Y} \widehat{Q}_n(Y)\widehat{\Phi}(Y)
\end{eqnarray}
where $\widehat{cap}$ is the capacity with respect to the weights $\widehat{Q}_n(Y)\widehat{r}_n(Y,Y')$ and $\widehat{\Phi}$ is the unique harmonic function with respect to the same weights with boundary values $\widehat{\Phi}|_{\pi_2^{-1}(m)}\equiv1$, $\widehat{\Phi}|_{\pi_2^{-1}(M)}\equiv0$. Due to this formula the sharp asymptotics of $\widehat{cap}(\pi_2^{-1}(m),\pi_2^{-1}(M))$ turns out to be the main ingredient in the sharp asymptotics of $\ev^{\pi_2^{-1}(m)}[\tau_{\pi_2^{-1}(M)}]$. We prove first an upper bound on $\widehat{cap}(\pi_2^{-1}(m),\pi_2^{-1}(M))$ in this subsection and show that it is sharp in the subsequent one.\\\quad\\ 
The upper bound involves the following quantities. Let $\widehat{\Lambda}'_n$ be the set of points in $\{Y:\;\widehat{I}(Y)<\infty\}$ which are at least of distance $n^{-\frac{1}{3}}$ from the boundary of the latter set. For each $Y\in\widehat{\Lambda}'_n$ let $D_n(Y)$ be the dice with center $Y$, edge length $2n^{-\frac{1}{2}}$ and edges parallel to an orthonormal basis of the eigenvectors of the Hessian $\widehat{H}(Y)$ of $\widehat{I}$ at $Y$. Moreover, define  
\begin{eqnarray}
\widehat{\scr Z}=\{z'^{(1)},\dots,z'^{(K')}\}=\pi_2^{-1}(\scr Z), 
\end{eqnarray}
let $\{z^{(1)},\dots,z^{(K)}\}$ be a minimal subset of $\widehat{\scr Z}$ with the property
\begin{eqnarray}
\widehat{\scr Z}\cap\widehat{\Lambda}'_n\subset\bigcup_{k=1}^K D_n(z^{(k)})
\end{eqnarray}
and denote $D_n(z^{(k)})$ by $D_n^{(k)}$. The idea behind this construction is that it will suffice to approximate the harmonic function on $\widehat{\Lambda}_n$ merely on $\bigcup_{k=1}^K D_n^{(k)}$ in order to obtain the sharp asymptotics of the capacity.\\\quad\\
Furthermore, define the eigenvalues $\gamma^{(k)}_i$, $1\leq i\leq L$ of the Hessian $\widehat{H}$ of $\widehat{I}$ at $z^{(k)}$ for $1\leq k\leq K$, the corresponding orthonormal bases of eigenvectors $v^{(k)}_i$, $1\leq i\leq L$, the possible transition directions $e_l$, $1\leq l\leq L$ at the saddle points, the corresponding transition probabilities $r^{(k)}_l$, the smallest (negative) eigenvalue $\lambda^{(k)}$ of the $L\times L$ matrix $\Big(\sqrt{r_l^{(k)}}\widehat{H}(z^{(k)})\sqrt{r_{l'}^{(k)}}\Big)_{l,l'}$ in the sense that the quadratic form given by $\widehat{H}(z^{(k)})$ and evaluated at pairs $(e^{(k)}_l,e^{(k)}_{l'})$ is rescaled by $\sqrt{r_l^{(k)}}\sqrt{r_{l'}^{(k)}}$, a corresponding normalized eigenvector $\widehat{w}^{(k)}$, $w^{(k)}\equiv\Big(\frac{\widehat{w}^{(k)}(l)}{\sqrt{r^{(k)}_l}}\Big)_l$ again in the sense that $\left\langle \widehat{w}^{(k)},e^{(k)}_l\right\rangle$ is rescaled by $\frac{1}{\sqrt{r^{(k)}_l}}$ and the sets
\begin{eqnarray}
\Gamma_{z^{(k)}}=\Big\{i:\;\gamma^{(k)}_i+2|\lambda^{(k)}|\left\langle w^{(k)},v^{(k)}_i\right\rangle^2\neq0\Big\}.
\end{eqnarray} 
With these notations we have  
{\lemma Under the assumptions of Theorem 1.2 it holds
\begin{eqnarray*}
\widehat{cap}(\pi_2^{-1}(m),\pi_2^{-1}(M))\leq\frac{n}{4\pi}
\sum_{k=1}^K \widehat{Q}_n(z^{(k)})|D^{(k)}_n\cap\widehat{\Lambda}_n|\sqrt{\frac{\pi}{2}}^{|\Gamma_{z^{(k)}}|}|\lambda^{(k)}|\\
\cdot\Big|\prod_{i\in \Gamma_{z^{(k)}}}\Big(\gamma^{(k)}_i+2|\lambda|\left\langle w^{(k)},v^{(k)}_i\right\rangle^2\Big)\Big|^{-\frac{1}{2}}
\sum_{l=1}^L r^{(k)}_l \left\langle e_l,v^{(k)}_1\right\rangle^2
(1+o(1))
\end{eqnarray*}
in the limit $n\rightarrow\infty$.}\\\quad\\
{\it Proof.} 1) As will become apparent from the proof, it suffices to treat the case $K=1$. In the proof of this case we will suppress the superscript $(1)$, i.e. write $\gamma_1$ for $\gamma^{(1)}_1$, $D_n$ for $D^{(1)}_n$ etc.\\\quad\\  
2) As was shown in \cite{begk3} it suffices to bound the capacity in a small neighborhood of $z$. We choose this set to be the dice $D_n$ as in the statement of the lemma with edge length modified to $2n^{-\frac{1}{2}+\eta}$ for a positive $\eta\ll 1$. With an abuse of notation we call the modified dice $D_n$ throughout the proof. Define the function
\begin{eqnarray}
g(Y)=\sqrt{\frac{n|\lambda|}{2\pi}}\int_{-\infty}^{\left\langle Y-z,w\right\rangle} e^{-n|\lambda|u^2/2}\;du
\end{eqnarray}
on $D_n$. We observe that in the limit $n\rightarrow\infty$ the function $g$ converges exponentially fast to $1$ on 
\begin{eqnarray*}
A_1=\{Y\in D_n:\;\left\langle Y-z,w\right\rangle>0\} 
\end{eqnarray*}
and to $0$ on
\begin{eqnarray*}
A_0=\{Y\in D_n:\;\left\langle Y-z,w\right\rangle<0\}.
\end{eqnarray*}
In particular, the value of $\widehat{cap}(\pi_2^{-1}(m),\pi_2^{-1}(M))$ is asymptotically upper bounded by the Dirichlet form restricted to $D_n$ and evaluated at $g$. Denoting inequalities up to a factor $1+o(1)$ by $\preceq$ and equalities up to the same factor by $\sim$ we have by the large deviation principle:
\begin{eqnarray*}
\widehat{cap}(\pi_2^{-1}(m),\pi_2^{-1}(M))\preceq\frac{1}{2}\sum_{Y,Y'\in D_n} \widehat{Q}_n(Y)\widehat{r}_n(Y,Y')\frac{n|\lambda|}{2\pi}
\left(\int_{\left\langle Y'-z,w\right\rangle}^{\left\langle Y-z,w\right\rangle}e^{-n|\lambda|u^2/2}\;du\right)^2\\
\sim\frac{n|\lambda|}{4\pi}\sum_{Y,Y'\in D_n} \widehat{Q}_n(z)\exp\left(-\frac{n}{2}\left\langle Y-z,\widehat{H}(z)(Y-z)\right\rangle\right)\widehat{r}_n(Y,Y')
\left(\int_{\left\langle Y'-z,w\right\rangle}^{\left\langle Y-z,w\right\rangle}e^{-n|\lambda|u^2/2}\;du\right)^2.
\end{eqnarray*}  
3) Next, we show the following two claims: $r_l$ is bounded away from zero for all $l$ with probability exponentially close to $1$ and for all $l$
\begin{eqnarray}
\sup_{Y\in D_n}\left|\frac{\widehat{r}_n(Y,Y+e_l)}{r_l}-1\right|\rightarrow_{n\rightarrow\infty}0.
\end{eqnarray}
To see this, assume without loss of generality that $e_l$ corresponds to a flip from $+$ to $-$, i.e. a decrease of $Y_a^+$ by $\frac{1}{n}$ and an increase of $Y_a^-$ by $\frac{1}{n}$ for some $a\in A$. We now claim that in
\begin{eqnarray}
r_l=\exp(-\beta n(v(\pi_2(z))-v(\pi_2(z+e_l)))_+)z_a^+
\end{eqnarray}
we have
\begin{eqnarray}
z_a^+=\frac{1}{n}|\{i:(\xi_1(i),\dots,\xi_p(i))=a,\sigma(i)=1\}|\geq C_1 
\end{eqnarray}
for a constant $C_1>0$ which would imply that $r_l$ is bounded away from $0$. This follows from 
\begin{eqnarray}
\frac{1}{n}|\{i:(\xi_1(i),\dots,\xi_p(i))=a\}|\geq C_2 
\end{eqnarray}
for a constant $C_2>0$ and the fact that $r_l$ has the same asymptotics as the transition probability at $z$ in the direction $-e_l$ due to reversibility of the dynamics and the second claim. Moreover, for any $Y\in D_n$ it holds 
\begin{eqnarray}
\widehat{r}_n(Y,Y+e_l)=\exp(-\beta n(v(\pi_2(Y))-v(\pi_2(Y+e_l)))_+)Y_a^+
\end{eqnarray}
and
\begin{eqnarray}
|Y_a^+-z_a^+|\leq C_3n^{-\frac{1}{2}+\eta}
\end{eqnarray}
with a constant $C_3>0$. This gives the second claim.\\\quad\\ 
4) Using this, the orthogonal decomposition of $\widehat{H}(z)$ and the leading term of the Taylor expansion of the integral, the asymptotics of the upper bound computes to
\begin{eqnarray*}
\frac{n|\lambda|}{4\pi}\sum_{Y\in D_n} \sum_{l=1}^L \widehat{Q}_n(z)\exp\left(-\frac{n}{2}\sum_{j=1}^L \gamma_j\left\langle Y-z,v_j\right\rangle^2\right)
r_l\exp\left(-n|\lambda|\left\langle Y-z,w\right\rangle^2\right)\left\langle e_l,w\right\rangle^2.
\end{eqnarray*}
By the standard approximation of integrals by Riemann sums the last expression is asymptotically equivalent to
\begin{eqnarray*}
\frac{n|\lambda|}{4\pi}\widehat{Q}_n(z)\sum_{l=1}^L r_l \left\langle e_l,w\right\rangle^2 \frac{|D_n\cap\widehat{\Lambda}_n|}{2^L n^{-\frac{L}{2}+L\eta}}\int_{D_n} 
\exp\left(-\frac{n}{2}\left(\sum_{j=1}^L \gamma_j\left\langle y,v_j\right\rangle^2+2|\lambda|\left\langle y,w\right\rangle^2\right)\right)\;dy\\
=\frac{n|\lambda|}{4\pi}\widehat{Q}_n(z)\sum_{l=1}^L r_l \left\langle e_l,w\right\rangle^2 \frac{|D_n\cap\widehat{\Lambda}_n|}{2^{|\Gamma_z|}n^{|\Gamma_z|(-\frac{1}{2}+\eta)}}\left(\sqrt{\frac{2\pi}{n}}\right)^{|\Gamma_z|}
\Big|\prod_{i\in \Gamma_z}\Big(\gamma_i+2|\lambda|\left\langle w,v_i\right\rangle^2\Big)\Big|^{-\frac{1}{2}}.
\end{eqnarray*}
Taking the limit $\eta\downarrow0$ we obtain the lemma in the case $K=1$. For the general case it suffices to note that the same construction can be applied for every $D_n^{(k)}$, $1\leq k\leq K$ and that the contribution of the set of points which belong to more than one dice is negligible in the limit.\ep\\\quad\\
{\bf 2.3. Sharpness of the upper bound.} Our goal is to show that the upper bound of Lemma 2.1 gives the exact asymptotics of the capacity. One of the main contributions of this paper compared to \cite{bbi} is that in our proof we do not construct an approximately harmonic flow to get a sharp lower bound from the Berman-Konsowa variational principle. Instead, we give a direct argument which shows that the function $g$ defined in the proof of lemma 2.1 approximates the unique harmonic function $\widehat{\Phi}$ on $\bigcup_{k} D^{(k)}_n$ with boundary values $1$ on $\bigcup_k A_1^{(k)}$ and $0$ on $\bigcup_k A_0^{(k)}$. To this end, we show in the next lemma how to derive an upper bound on $\|g-\widehat{\Phi}\|_{\infty}$ from an upper bound on $\|\Delta(g-\widehat{\Phi})\|_{\infty}$. In the proof of Theorem 2.3 we show that this bound is good enough to conclude that the Dirichlet form evaluated at $g$ gives the leading term of $\widehat{cap}(\pi_2^{-1}(m),\pi_2^{-1}(M))$. 
{\lemma Let $F$ be a function on $B\subset\widehat{\Lambda}_n$ and $B_0\subset B$ be such that each connected component of $B$ contains at least one point of $B_0$. Moreover, let the diameter with respect to the graph distance of each connected component of $B$ be bounded by $C$ for some $C>0$,
\begin{eqnarray*}
\|\Delta F\|_{\infty}\leq\delta,\quad |F|\leq\epsilon\;{\rm on}\;B_0   
\end{eqnarray*}
for some $\delta,\epsilon>0$ and $\Delta$ defined with respect to some weights on the edges which are uniformly bounded below by $\theta$. Then 
\begin{eqnarray}
\|F\|_{\infty}\leq \epsilon+\frac{C\delta}{\theta}.
\end{eqnarray}}
{\it Proof.} Without loss of generality we may assume that $B$ is connected (otherwise we prove the bound for each connected component separately). Moreover, it suffices to show $\max F\leq \epsilon+\frac{C\delta}{\theta}$, since the same argument can be then applied to $-F$. For $1\leq i\leq C$ define disjoint sets $B_i\subset B$ recursively by letting $B_i$ be the set of all points in $B$ which are connected to at least one point in $B_{i-1}$ and do not belong to $B_1,\dots,B_{i-1}$. Let $s_i$ be the maximum of $F$ on $B_i$ and $S$ be the maximum of $F$ on $B$. Moreover, let $j$ be such that $s_j=S$ and assume without loss of generality $j\neq0$. Then there exists a $Y_j\in B_j$ with $F(Y_j)=S$. Using the bound on $|(\Delta F)(Y_j)|$, $Y_j\in B_j$ and the fact that the weights are bounded below by $\theta$, we conclude that there exists a $Y_{j-1}\in B_{j-1}$ with $F(Y_{j-1})\geq S-\frac{\delta}{\theta}$. Denoting the (not necessarily unique) $\arg\max$ of $F$ on $B_{j-1}$ by $Y_{j-1}$ we have $s_{j-1}=F(Y_{j-1})\geq S-\frac{\delta}{\theta}$. Iterating the argument $j$ times we end up with $s_0\geq S-j\frac{\delta}{\theta}\geq S-C\frac{\delta}{\theta}$. On the other hand, $s_0\leq\epsilon$. Thus, $S\leq \epsilon+\frac{C\delta}{\theta}$ as claimed.\ep\\\quad\\
Equipped with the lemma we are ready to prove
{\theorem Under the assumptions of Theorem 1.2 it holds
\begin{eqnarray*}
\widehat{cap}(\pi_2^{-1}(m),\pi_2^{-1}(M))=\frac{n}{4\pi}
\sum_{k=1}^K \widehat{Q}_n(z^{(k)})|D^{(k)}_n\cap\widehat{\Lambda}_n|\sqrt{\frac{\pi}{2}}^{|\Gamma_{z^{(k)}}|}|\lambda^{(k)}|\\
\cdot\Big|\prod_{i\in \Gamma_{z^{(k)}}}\Big(\gamma^{(k)}_i+2|\lambda|\left\langle w^{(k)},v^{(k)}_i\right\rangle^2\Big)\Big|^{-\frac{1}{2}}
\sum_{l=1}^L r^{(k)}_l \left\langle e_l,v^{(k)}_1\right\rangle^2
(1+o(1))
\end{eqnarray*}
in the limit $n\rightarrow\infty$.}\\\quad\\
{\it Proof.} 1) It will become apparent from the proof that we may restrict to the case $K=1$ and as in the proof of Lemma 2.1 we will again suppress the superscript $(1)$. Moreover, for $Y\in D_n$ we define 
\begin{eqnarray}
\widehat{Q}'_n(Y)=\widehat{Q}_n(z)\exp\Big(-\frac{n}{2}\left\langle Y-z,\widehat{H}(z)(Y-z)\right\rangle\Big)
\end{eqnarray}
and notice
\begin{eqnarray}
\Big|\frac{\widehat{Q}_n(Y)}{\widehat{Q}'_n(Y)}-1\Big|\leq\frac{C_1}{\sqrt{n}}
\end{eqnarray}
for a constant $C_1>0$. Thus, step 3 of the proof of lemma 2.1 together with lemma 4.1 of \cite{bbi} imply
\begin{eqnarray}
\widehat{cap}(\pi_2^{-1}(m),\pi_2^{-1}(M))=\widehat{cap}'(\pi_2^{-1}(m),\pi_2^{-1}(M))(1+o(1))
\end{eqnarray}
where $\widehat{cap}'$ is defined with respect to the weights $\widehat{Q}'_n(Y)r_l$ on $D_n$. Finally, let $\widehat{d}'$ be the Dirichlet form on $D_n$ with respect to the same weights and $\widehat{\Phi}$ be the restriction to $D_n$ of the harmonic function on the enlarged dice with edge length $2n^{-\frac{1}{2}+\eta}$ with boundary values $0$ on $A_0$ and $1$ on $A_1$ where $\eta\ll 1$ is fixed and $A_0$, $A_1$ are defined as in the proof of Lemma 2.1.\\\quad\\
2) Let $e'_l$ be the unit vector in the same direction as $e_l$. As in lemma 4.4 of \cite{bbi} one shows that
\begin{eqnarray}
|\Delta(g-\widehat{\Phi})(Y)|\leq\widehat{Q}'_n(Y)\sqrt{\frac{|\lambda|}{2\pi n}}e^{-n|\lambda|\left\langle Y-z,w\right\rangle^2/2}
\Big|\sum_{l=1}^L r_l\left\langle e'_l,w\right\rangle\Big|\frac{C_2}{n}\equiv\delta
\end{eqnarray}
for a constant $C_2>0$ where $\Delta$ is defined with respect to the weights $\widehat{Q}'_n(Z)r_l$. Next, let
\begin{eqnarray}
\epsilon\equiv\sup_{\partial D_n} |g-\widehat{\Phi}|,\\
\theta\equiv\min_{Y\in D_n,1\leq l\leq L} \widehat{Q}'_n(Y)r_l.
\end{eqnarray}
We may now apply lemma 2.2 for $F=g-\widehat{\Phi}$ with $B=D_n$, $\delta$, $\epsilon$, $\theta$ as just defined and noting that $C=O(\sqrt{n})$ in our case. This gives
\begin{eqnarray}
|g-\widehat{\Phi}|_{\infty}\leq \epsilon+C_3 \sup_{Y\in D_n} \frac{\sqrt{n}}{\theta}\widehat{Q}'_n(Y)\left(\sqrt{\frac{|\lambda|}{2\pi n}}
e^{-n|\lambda|\left\langle Y-z,w\right\rangle^2/2}\Big|\sum_{l=1}^L r_l\left\langle e'_l,w\right\rangle\Big|\right)\frac{1}{n}
\end{eqnarray}
for a constant $C_3>0$. Since the edge length of $D_n$ is $2n^{-\frac{1}{2}}$, we deduce
\begin{eqnarray*}
\sup_{Y\in D_n}\frac{1}{\theta}\widehat{Q}'_n(Y)e^{-n|\lambda|\left\langle Y-z,w\right\rangle^2/2}
\leq\frac{1}{\min_{1\leq l\leq L} r_l}\sup_{Y,Y'\in D_n}\frac{\exp\Big(-\frac{n}{2}\left\langle Y-z,\widehat{H}(z)(Y-z)\right\rangle\Big)}
{\exp\Big(-\frac{n}{2}\left\langle Y'-z,\widehat{H}(z)(Y'-z)\right\rangle\Big)}\\
\leq\frac{1}{\min_{1\leq l\leq L} r_l}\cdot e^{\frac{L}{2}(\gamma_L-\gamma_1)}.
\end{eqnarray*}
Moreover, we know from the proof of lemma 2.1 that $\min_{1\leq l\leq L} r_l$ is bounded away from $0$, so the latter two bounds imply
\begin{eqnarray}
|g-\widehat{\Phi}|_{\infty}\leq\frac{C_4}{n}
\end{eqnarray}
for a $C_4>0$.\\\quad\\
3) By the uniform bound on $|g-\widehat{\Phi}|$ we obtain for all $Y\in D_n$ and all $1\leq l\leq L$:
\begin{eqnarray*}
|(g(Y+e_l)-g(Y))^2-(\widehat{\Phi}(Y+e_l)-\widehat{\Phi}(Y))^2|\leq\frac{2C_4}{n}|g(Y+e_l)-g(Y)+\widehat{\Phi}(Y+e_l)-\widehat{\Phi}(Y)|\\
\leq\frac{2C_4}{n}\Big(\frac{2C_4}{n}+2|g(Y+e_l)-g(Y)|\Big)
=\frac{2C_4}{n}\left(\frac{2C_4}{n}+\sqrt{\frac{2}{\pi}}\int_{\sqrt{n|\lambda|}\left\langle Y+e_l-z,w\right\rangle}^{\sqrt{n|\lambda|}\left\langle Y-z,w\right\rangle}e^{-v^2/2}\;dv\right)
\leq C_5 n^{-\frac{3}{2}}
\end{eqnarray*}
with a constant $C_5>0$, since the length of the interval of integration is of order $\frac{1}{\sqrt{n}}$. This implies directly
\begin{eqnarray}
\frac{1}{\widehat{d}'(g)}|\widehat{d}'(g)-\widehat{d}'(\widehat{\Phi})|
\leq C_5 n^{-\frac{3}{2}}\frac{\sum_{Y\in D_n,1\leq l\leq L}\widehat{Q}'_n(Y)r_l}{2\widehat{d}'(g)}.
\end{eqnarray}
Finally, we claim 
\begin{eqnarray}
\frac{\sum_{Y\in D_n, 1\leq l\leq L}\widehat{Q}'_n(Y)r_l}{2\widehat{d}'(g)}\leq C_6 n  
\end{eqnarray}
for a $C_6>0$. To this end, we observe that in the denominator we are summing the terms
\begin{eqnarray*}
\widehat{Q}'_n(Y)r_l\frac{n|\lambda|}{2\pi}
\left(\int_{\left\langle Y+e_l-z,w\right\rangle}^{\left\langle Y-z,w\right\rangle}e^{- n|\lambda|u^2/2}\;du\right)^2
=\widehat{Q}'_n(Y)r_l\frac{1}{2\pi}\left(\int_{\sqrt{n|\lambda|}
\left\langle Y+e_l-z,w\right\rangle}^{\sqrt{n|\lambda|}\left\langle Y-z,w\right \rangle} e^{-v^2/2}\;dv\right)^2
\end{eqnarray*}
over the same set as in the numerator. Since $\sqrt{n|\gamma_1|}|\left\langle Y-z,v_1\right\rangle|$ is bounded and the length of the interval of integration in the latter integral is of order $\frac{1}{\sqrt{n}}$, the integral itself is of order $\frac{1}{\sqrt{n}}$, so the fraction is of order $n$ as claimed. All in all, we have shown that there exists a $C_7>0$ such that
\begin{eqnarray}
\frac{1}{\widehat{d}'(g)}|\widehat{d}'(g)-\widehat{d}'(\widehat{\Phi})|\leq \frac{C_7}{\sqrt{n}}
\end{eqnarray}
which together with $\widehat{cap}'(\pi_2^{-1}(m),\pi_2^{-1}(M))=\widehat{d}'(\widehat{\Phi})$ and (II.33) implies that the leading term of $\widehat{cap}(\pi_2^{-1}(m),\pi_2^{-1}(M))$ is given by $\widehat{d}'(g)$. Noting that the upper bound of lemma 2.1 is precisely $\widehat{d}'(g)$ we deduce the theorem.\ep 

\section{Proof of the main result}

In this last section we prove the following Proposition 3.1 which yields directly Theorem 1.2 and provides the value of the correction term $c(n)$ therein.
{\prop If assumption 1.1 is satisfied, then Theorem 1.2 holds with
\begin{eqnarray*}
c(n)=\frac{4\pi}{n}\sqrt{\frac{\pi}{2}}^{|\Gamma_{m^{(1)}}|}\Big|\prod_{i\in\Gamma_{m^{(1)}}}\widetilde{\gamma}_i^{(1)}\Big|^{-\frac{1}{2}}\widehat{\kappa}_n(m^{(1)})\\
\cdot\left(\sum_{k=1}^K \widehat{\kappa}_n(z^{(k)})\sqrt{\frac{\pi}{2}}^{|\Gamma_{z^{(k)}}|}|\lambda^{(k)}|
\Big|\prod_{i\in \Gamma_{z^{(k)}}}\Big(\gamma^{(k)}_i+2|\lambda|\left\langle w^{(k)},v^{(k)}_i\right\rangle^2\Big)\Big|^{-\frac{1}{2}}
\sum_{l=1}^L r^{(k)}_l \left\langle e_l,v^{(k)}_1\right\rangle^2\right)^{-1}
\end{eqnarray*}
where $m^{(1)}$ is an element of $\pi_2^{-1}(m)$, $\widetilde{\gamma}_i^{(1)}$, $1\leq i\leq L$ are the eigenvalues of $\widehat{H}(m^{(1)})$ and
$\Gamma_{m^{(1)}}=\{i:\widetilde{\gamma}_i^{(1)}\neq0\}$.}\\\quad\\ 
{\it Proof.} 1) Having established the formula
\begin{eqnarray}
\ev^{\pi^{-1}(m)}[\tau_{\pi^{-1}(M)}]=\ev^{\pi_2^{-1}(m)}[\tau_{\pi_2^{-1}(M)}]=\frac{1}{\widehat{cap}(\pi_2^{-1}(m),\pi_2^{-1}(M))}\sum_{Y} \widehat{Q}_n(Y)\widehat{\Phi}(Y)
\end{eqnarray}
and Theorem 2.3, it remains to compute the sharp asymptotics of $\sum_{Y} \widehat{Q}_n(Y)\widehat{\Phi}(Y)$. To this end, for $a>0$ define $U^a(m)$ to be the union of connected components of
\begin{eqnarray*}
\{Y\in\widehat{\Lambda}_n:\;\widehat{I}(Y)\leq I(m)+a\}
\end{eqnarray*}
which have a non-empty intersection with $\pi_2^{-1}(m)$. Moreover, let $U^a(M)$ be the rest of this set. Choosing $a$ small enough, the two sets become disjoint neighborhoods of $\pi_2^{-1}(m)$ and $\pi_2^{-1}(M)$, respectively. Following section 6 of \cite{bbi} one shows that in the limit $n\rightarrow\infty$ the only relevant contribution of $\sum_{Y} \widehat{Q}_n(Y)\widehat{\Phi}(Y)$ comes from $U^a(m)$ and on this set the summands can be replaced by $\widehat{Q}_n(Y)$ without changing the asymptotics.\\\quad\\ 
2) For the sake of completeness we give an outline of the argument of the first claim here. The second claim is completely analogous. Since $\widehat{Q}_n$ is small outside $U^{a}(m)\cup U^a(M)$, it suffices to show that $\widehat{\Phi}$ is small on $U^{a}(M)$. To this end let $A=\pi_2^{-1}(m)$ and $B=\pi_2^{-1}(M)$. Note that for any $Y\in\widehat{\Lambda}_n-(A\cup B)$ and any set $C\subset\widehat{\Lambda}_n-(A\cup B)$ the Markov property imlies
\begin{eqnarray}
\widehat{\Phi}(Y)=\pp^Y(\tau_A<\tau_B)\leq\pp^Y(\tau_A<\tau_{B\cup C})+\pp^Y(\tau_C<\tau_B)\cdot\max_{Y'\in C}\pp^{Y'}(\tau_A<\tau_B). 
\end{eqnarray}
As observed in Proposition 6.12 of \cite{bbi} this implies the claim provided we can find a set $B\subset D\subset U^a(M)$ and, choosing $C$ to be the part of the boundary of $D$ which does not belong to $B$, find constants $0<c_1<c_2$, $0<c_3$ with
\begin{eqnarray}
\forall Y\in C:\quad \pp^Y(\tau_C<\tau_B)\leq 1-e^{-c_1 n},\\
\forall Y\in C:\quad \pp^Y(\tau_A<\tau_{B\cup C})\leq e^{-c_2 n},\\
\forall Y\in U^{a}(M)-D:\quad \pp^Y(\tau_A<\tau_{B\cup C})\leq e^{-c_3 n}. 
\end{eqnarray} 
The main part of the proof of the claim consists of proving the third inequality. In fact, the other two can be deduced from it using bounds on the derivatives of $I$ which in our case can be replaced by positive constants, since in our case we can find a neighborhood of the critical points of $\widehat{I}$ contained in a compact set
\begin{eqnarray*}
K\subset\{Y:\;\widehat{I}(Y)<\infty\}
\end{eqnarray*}
which is independent of $n$. The proof of the third inequality consists of showing that on $U^{a}(M)-D$
\begin{eqnarray}
\pp^Y(\tau_A<\tau_{B\cup C})\leq\max_{Y'\in C}\frac{\Psi(\sigma)}{\Psi(\sigma')}
\end{eqnarray}  
for any super-harmonic function $\Psi$ on $U^{a}(M)-D$ (see p. 45 in \cite{bbi}) and that for any $\alpha\in(0,1)$
\begin{eqnarray}
\Psi(Y)=e^{(1-\alpha)n I(\pi_2(Y))}
\end{eqnarray}  
is super-harmonic for large $n$. The latter can be done as in Proposition 6.4 of \cite{bbi}. Hereby, we do not need to treat the case that $\pi_2(Y)$ is in a small neighborhhood of the boundary of $\{x\in[-1,1]^p:\;I(x)<\infty\}$ separately, since $U^a(M)\subset K$ for $a$ small enough.\\\quad\\
3) Finally, one may replace $U^a(m)$ by a union of dices $D^b(m^{(s)})$, $1\leq s\leq S$ around 
\begin{eqnarray}
\{m^{(1)},\dots,m^{(S)}\}=\pi_2^{-1}(m)
\end{eqnarray}
with edge length $b\ll a$ and edges parallel to an orthonormal basis of eigenvectors of $\widehat{H}(m^{(s)})$, since such a union contains $U^{\widetilde{a}}(m)$ for an $\widetilde{a}<a$ small enough. Thus, by approximating the sum with an integral and evaluating the integral as in the proof of lemma 2.1 we obtain     
\begin{eqnarray*}
\sum_Y\widehat{Q}_n(Y)\widehat{\Phi}(Y)\sim \sum_{s=1}^S \widehat{Q}_n(m^{(s)})\Big(\sum_{Y\in D^b(m^{(s)})}\exp\left(-\frac{n}{2}\left\langle Y-m^{(s)},\widehat{H}(m^{(s)})(Y-m^{(s)})\right\rangle\right)\Big)\\
\sim\sum_{s=1}^S \widehat{Q}_n(m^{(s)})|D^b(m^{(s)})\cap\widehat{\Lambda}_n|b^{-|\Gamma_{m^{(s)}}|}(2\pi)^{\frac{|\Gamma_{m^{(s)}}|}{2}}
n^{-\frac{|\Gamma_{m^{(s)}}|}{2}}\Big|\prod_{i\in\Gamma_{m^{(s)}}}\widetilde{\gamma}_i^{(s)}\Big|^{-\frac{1}{2}}
\end{eqnarray*}   
where $\widetilde{\gamma}_i^{(s)}$, $1\leq i\leq L$ are the eigenvalues of $\widehat{H}(m^{(s)})$ for $1\leq s\leq S$ and
\begin{eqnarray}
\Gamma_{m^{(s)}}=\{i:\;\widetilde{\gamma}_i^{(s)}\neq0\}.
\end{eqnarray}
We would like to take the limit $b\downarrow0$ of the latter expression. This cannot be done directly, since the evaluation of the integral in the last step would become invalid. However, we can take the limit along a sequence $b(n)=2n^{-\frac{1}{2}+\eta}$ for a positive $\eta\ll 1$ and take the limit $\eta\downarrow0$ in the end. By this means, we obtain
\begin{eqnarray*}
\sum_{s=1}^S \widehat{Q}_n(m^{(s)})|D^{b(n)}(m^{(s)})\cap\widehat{\Lambda}_n|\sqrt{\frac{\pi}{2}}^{|\Gamma_{m^{(s)}}|}
\Big|\prod_{i\in\Gamma_{m^{(s)}}}\widetilde{\gamma}_i^{(s)}\Big|^{-\frac{1}{2}}
\end{eqnarray*}
as the leading term of $\sum_Y\widehat{Q}_n(Y)\widehat{\Phi}(Y)$. Since the summands do not depend on $s$, the latter sum evaluates to
\begin{eqnarray*}
Q_n(m)|D^{b(n)}(m^{(1)})\cap\widehat{\Lambda}_n|\sqrt{\frac{\pi}{2}}^{|\Gamma_{m^{(1)}}|}
\Big|\prod_{i\in\Gamma_{m^{(1)}}}\widetilde{\gamma}_i^{(1)}\Big|^{-\frac{1}{2}}.
\end{eqnarray*}
Putting this together with Theorem 2.3 and equation (III.31) gives Proposition 3.1 by noting that $|D^{b(n)}(m^{(1)})\cap\widehat{\Lambda}_n|$ has the same asymptotics as $|D^{(k)}_n\cap\widehat{\Lambda}_n|$ for all $1\leq k\leq K$ and by using the large deviation principle for $Q_n$.\ep \\\quad\\\quad\\
{\bf Acknowledgement.} The author thanks Amir Dembo for many helpful discussions and suggestions throughout the preparation of this work.

\end{document}